\documentclass[a4paper,11pt]{article}

\usepackage{amsmath,english,german,latexsym,amssymb,curves}
\newtheorem{thm}{Theorem}
\newtheorem{prop}[thm]{Proposition}
\newtheorem{lemma}{Lemma}
\newtheorem{cor}{Corollary}
\newtheorem{main}{Main Theorem}

\newcommand{\mtiny}{\scriptscriptstyle}

\newcommand{\proof}[1][]{{\it Proof#1: }}

\newcommand{\qed}[1][1cm]{\hspace*{\fill} $\Box$ \vspace{#1}}

\renewcommand{\P}{{{\mathbf P}^1}}

\newcommand{\CC}{{\mathbf C}}

\newcommand{\del}{\partial}

\newcommand{\id}{\operatorname{id}}

\renewcommand{\a}{\alpha}
\renewcommand{\b}{\beta}
\newcommand{\cg}{\gamma}

\renewcommand{\d}{\delta}
\newcommand{\D}{\Delta}

\renewcommand{\l}{\lambda}

\newcommand{\s}{\sigma}

\newcommand{\tto}{\longrightarrow}

\newcommand{\Feins}{{\bf F}_{1}}
\newcommand{\Fzwei}{{\bf F}_{2}}

\newcommand{\fami}{{\cal F}}

\newcommand{\milnor}{{\cal U}}

\newcommand{\afami}{{\cal A}}
\newcommand{\cfami}{{\cal C}}
\newcommand{\efami}{{\cal E}}

\newcommand{\zfami}{{\cal Z}}

\newcommand{\yfami}{{\cal Y}}

\newcommand{\labell}[1]{\label{#1}}

\newcommand{\G}{\Gamma}

\newcommand{\inv}{^{-1}}

\newcommand{\cq}{{\chi,q}}
\newcommand{\cqm}{_{\chi,q}^{\mu}}
\newcommand{\Xhalf}{X{\mtiny 1_{\!\!\;\!/_{\!2}}}}
\newcommand{\Fchi}{{\bf F}_{2\chi}}
\newcommand{\Gcq}{\Gamma_\cq}

\begin{document}

\begin{center}
{\Large\sc
Monodromy groups of irregular elliptic surfaces}\\[1.5cm]

Michael L\"onne
\footnote[1]{
M. L"onne: Institut f"ur Mathematik, Universit"at Hannover, Am Welfengarten
1,\\ 30167 Hannover, Deutschland. e-mail: loenne@math.uni-hannover.de}
\end{center}

\begin{abstract}

Monodromy groups, i.e.\ the groups of isometries of the intersection lattice
$L_X:=H_2/\!\:\!$\raisebox{-0.7mm}{torsion} generated by the monodromy action of
all deformation families of a given surface, have been computed in \cite{ellmon}
for any minimal elliptic surface with $p_g>q=0$. New and refined methods are now
employed to address the cases of minimal elliptic surfaces with
$p_g\geq q>0$.\\
Thereby we get explicit families such that any isometry is
in the group generated by their monodromies or does not respect the invariance of
the canonical class or the spinor norm.
The monodromy is also shown to act by the full symplectic group on the first
homology modulo torsion.

\end{abstract}

\section*
{Introduction}

Monodromy is a powerful tool and has been a predominant subject of interest in
the realms of singularity since long.
This paper as its predecessor tries to broaden the view to include the case of
deformation families of compact complex surfaces.
This first step leaving the ground of surfaces singularities and complete
intersection surfaces just reaches the case of elliptic surfaces but nevertheless
it establishes an astonishing likeness between the topological notion of the
representation of the diffeomorphism group and the analytical notion of the
monodromy action. In fact the associated isometry groups are the same up to index
two.\\
In the first section we construct families with monodromy acting on a simply
connected part of a reference surface which is isomorphic to the Milnor fibre of a
hypersurface singularity. But in contrast to the case of regular elliptic surfaces
this Milnor fibre supports only classes generating a sublattice of $L_X$ of
possibly high corank.\\
So next we give families which yield vanishing cycles generating a sublattice
containing the classes of Lagrangian tori in the bundle part of the reference
surface.\\
Finally we have to take care again of the fibre classes which is done as in the
regular case. The results thus collected fit together nicely and yield the main
result:

\begin{main}
\labell{main}
Let $X$ be a minimal elliptic surface with positive holomorphic Euler number
$\chi=p_g-q+1$ and positive irregularity $q$, then there exist families of elliptic
surfaces containing $X$, such that the induced monodromy actions on the homology
lattice
$L_X$ generate $O'_k(L_X)$, the group of isometries of real spinor norm one
fixing the canonical class.
\end{main}

%\newpage

\section*
{Families obtained from unfoldings}

The aim of the first section is to construct families of irregular elliptic
surfaces, which contain Morsifications of isolated singularities of maximal
possible Milnor numbers.\\
Let us start with the Hirzebruch surfaces $\Feins,\Fzwei$ ruled over $\P$
containing sections $\s^\infty_1,\s^\infty_2$ of negative square and with branch
divisors $B_1,B_2$ resp. given by 
\begin{eqnarray*}
B_1 & = & \s^\infty_1\cup \overline{\{x(y^3-3y-x^3)=0\}},\\
B_2 & = & \s^\infty_1\cup \overline{\{y^3-3y-x^6=0\}}.
\end{eqnarray*}
The corresponding double covers are rational elliptic surfaces which shall be
denoted by $\Xhalf,X_1$ respectively.

Generalizing this construction we define $\zfami_\chi$ to
be the family of elliptic surfaces given by the family of double covers of $\Fchi$
along the family of branch curves consisting of
the curves given by $y^3-3y-x^{6\chi}-\l x-u$ together with the negative section:
\begin{eqnarray*}
\zfami\hspace{4mm} & \to & \hspace{1mm}\P\times\CC^2\\
\searrow & & \swarrow\\
& \CC^2 & \ni(\l,u)
\end{eqnarray*}
The negative section lifts to $\s_\zfami$ which restricts to a section of each
surface of $\zfami$.

Note that incidentally for $\chi>1$ the surface $Z_0$ in $\zfami_\chi$ with
parameters
$\l=u=0$ coincides with the pull back of $X_1$ along the degree $\chi$ cyclic
Galois cover of $\P$ branched at zero and infinity. Similarly if we take the pull
back of
$\Xhalf$ along the double cover of $\P$ branched at zero and infinity we have to
normalize and blow down exceptional divisors in the fibres to get the central
surface $Z_0$ of $\zfami_1$.

\begin{lemma}
\labell{Sbase}
There is an open polycylinder $S\subset \CC^2$ containing the origin and an open
holomorphic disc $ D$ with $\infty\in D\subset\bar D\subset\P\setminus\{0\}$
such that
\begin{enumerate}
\item
for all $s\in S$ the branch points of the Kodaira $j$-function $j_s$ on $\P$
associated to $Z_s$ are contained in $ D\cup\{0\}$,
\item
for all $s\in S$ the preimage $U_s:=j_s\inv(D)$ is a holomorphic disc which
covers $D$ with degree $12\chi$,
\item
the induced family $\zfami'$ of elliptically fibred surfaces over varying
holomorphic discs
\begin{eqnarray*}
\zfami'\hspace{4mm} & \to & \hspace{4mm}\{U_s\}_{s\in S}\subset\P\times S\\
\searrow & & \swarrow\\
& S &
\end{eqnarray*}
contains a Morsification of a singularity of type $J_{2\chi}$.
\end{enumerate}
\end{lemma}

\proof
The $j$-function can explicitly be given in terms of the base coordinate $x$ and
the parameter $s=(u,\l)$ as
$$
j_s(x)=\frac1{1-4(x^{6\chi}+\l x+u)^2}.
$$
The formula shows that the associated maps from $\P$ to itself are of degree
$12\chi$ with no zero except at infinity. Hence zero is a branch point of constant
multiplicity and does not vary.\\
On the other hand the branching in the neighbourhood of infinity does vary, but
since branch points vary continuously with the parameters, their range is bounded
as soon as the parameters are.
Thus for any bounded polydisc $S$ a disc $D$ can be found as claimed such that
$i)$ holds.\\
To get $ii)$ we argue with the disc $D^c$ complementary to $D$ in $\P$. Since
$j_s$ is totally ramified at zero with no other branching in $D^c$ it is
equivalent to the standard branching $z\mapsto z^{12\chi}$, and the preimage is
therefore a disc. And so is its complement, the cover $U_s$ of $D$.\\
Finally we observe that in the family
$f_\l=z^2+y^3-3y-x^{6\chi}-\l x$ Morse functions are dense, cf.\ \cite{sing}, and that
any bound on
$\l$ will imply a uniform bound on the range of the critical values of the
$f_\l$.
Choose $S$ accordingly to be any polycylinder given by such a pair of bounds.
Then by construction $\zfami|_S$ contains a Morsification of the singularity of
type $J_{2\chi}$ given by the local equation $z^2+y^3=x^{6\chi}$.
Since each surface in $\zfami|_S$ fibres smoothly over $\P\setminus U_s$, also
$\zfami'$ contains the Morsification and so $iii)$ holds.
\qed

The construction of a reference surface uses suitable branched covers:

\begin{lemma}
\labell{hurwitz}
Given positive integers $n\!\geq\!2,q$ and a set $\{x_0,...,x_{2q+n}\}$ of
distinct points of $\P$,
then there is a curve $C_q$ of genus $q$ and a branched covering $\pi:C_q\to\P$
of degree $n$ totally ramified at $x_1$ such that the associated monodromy
representation
$$
\rho:\pi_1(\P\setminus\cup_{i>0}\{x_i\},x_0)\to{\bf S}_n
$$
assigns the following permutations to simple loops $\omega_i$ associated to a given
geometric basis of paths from $x_0$ to the $x_i$:
$$
\begin{array}{cccl}
\rho([\omega_j]) & = & (12) & j=2,...,2q+2\\
\rho([\omega_{2q+j}]) & = & (j\!-\!\!1\,\,\,j) & j=3,...,n.
\end{array}
$$
\end{lemma}

\proof
The given monodromy representation assigns to the composition of the loops
$\omega_i$, $i${}$>$1 an element of order $n$. Hence the unramified cover of
$\P\setminus\{x_1\}$ determined by the subgroup $\rho\inv({\bf S}_{n-1})$ of $\pi_1$
can be compactified by a single point to yield a branched covering $\pi:C_q\to\P$
as claimed.
\qed

\begin{prop}
\labell{e-ref}
Given discrete invariants $\chi,q,\mu=m_1,...,m_n$ there is a smooth elliptic
surface $X$ with the given invariants which contains the central surface $Z'_0$
of $\zfami'_\chi\to S$.
\end{prop}

\proof
Let us first consider the case $\chi>1$:
Given any simply connected neighbourhood $U$ of $x_1=0$ in $\P$ there is as in
lemma \ref{hurwitz} a cover $\pi:C_q\to\P$ of degree $\chi$ and genus $q$ totally
branched at zero with the additional property that all other branch points are in
the complement of $U$. Then $U$ is covered by a disc $\tilde U$ in $C_q$ with total
ramification at zero, just as in the restriction of the cyclic cover $c:\P\to\P$ of
degree $\chi$ branched at $\{0,\infty\}$ into which $\tilde U\to U$ thus
necessarily embeds.\\
If we let $X'$ be the pull back of $X_1$ along $\pi$ and
choose $U$ to contain
$c(U_0)$ then $Z'_0$ embeds holomorphically into $X'|_{\tilde U}$ as follows
from the observation made above.
This property is not affected by the logarithmic transformations on $X'$ necessary
to get $X$ if they are performed over the complement of $\tilde U$.\\
In the case of $\chi=1$ we let $X'$ instead be the minimal model of the
normalization of the pull back of $\Xhalf$ along a hyperelliptic cover $\pi$ with
branching outside $U$ except for the branch point at zero.
\qed

\begin{lemma}
\labell{unrami}
Let $Z_s$ and $X$ be elliptic surfaces as introduced above, then there is an
annulus $A$ obtained from $D$ as given by
lemma \ref{Sbase} by removing a closed disc such that
\begin{enumerate}
\item
$A_s:=j_s\inv(A)$ is an annulus for all $s\in S$ where $j_s$ is the $j$-invariant
on $\P$ associated to the elliptic surface $Z_s$,
\item
$A_X:=j_X\inv(A)$ is an annulus for the $j$-invariant $j_X$
on $C_q$ associated to the elliptic surface $X$.
\end{enumerate}
\end{lemma}

\proof
If $A$ is chosen suitably it contains no branchpoints for all $j_s$ and $j_X$
and thus is unramified covered by $A_s$ resp.\ $A_X$.
The coverings are connected since the boundary coverings $\del U_s\to\del D$ are
connected by lemma \ref{Sbase} and
$\del U_X\to\del D$, where $U_X=j_X\inv(D)$, is equivalent to
$\del U_0\to\del D$ by construction of $X$, hence the claim.
\qed

\begin{lemma}
\labell{annulus}
There is an embedding $h=(h_1,h_2)$ of the subset $A_S:=\cup_sA_s$ of $\P\times S$
into
$C_q\times S$ such that
$$
j_X(h_1(a,s))=j_s(a)=:j_S(a,s),\quad\forall (a,s)\in A_S.
$$
\end{lemma}

\proof
Both $A_S$ and $A_X\times S$ are connected unramified covers of $A\times S$ of
degree $12\chi$ via $j_S$ resp. $j_X\times\id_S$. Since the fundamental
group of
$A\times S$ is cyclic both covers are homotopically equivalent and therefore
holomorphically equivalent. The embedding is thus obtained by an equivalence
followed by the obvious embedding of $A_X\times S$ into $C_q\times S$.
\qed 

\begin{lemma}
\labell{iso-lift}
There is a unique holomorphic embedding $\tilde h$ of $\zfami|_{A_S}$
into $X\times S$ covering $h$
which maps points of
$\s_\zfami$ into $\s_0\times S\subset Z'_0\times S\subset X\times S$.
\end{lemma}

\proof
The functional invariant associated to $\del U_s$ varies continuously with $s\in S$
and thus is constant. The functional invariant along $\del U_X$ is the same, for
it is equal to that of $\del U_0$ by construction of $X$.
Hence the functional invariant and the \mbox{$j$-invariant} by lemma \ref{annulus}
coincide on
$A_S$ and $A_X\times S$, and we get with thm.\ I.3.19 of \cite{fm} a unique
embedding $\tilde h$ of the fibrations covering $h$ which respects the sections.
\qed

By discarding a suitable neighbourhood of the singular fibres in $X\times S$, the
map $\tilde h$ becomes the biholomorphic identification of collars. By gluing we
thus get a family $\fami\cqm\to S$ of elliptic surfaces.\\
Notice that $X$ is biholomorphic to the central surface and let $\milnor$ denote the
open subset of $\fami$ given by $\zfami'\setminus\s_\zfami$.\\

Our final aim in this paragraph is to describe the algebraic monodromy group
$\G(\fami)\subset O(L_{X})$ of the family $\fami$
in terms of a set $\D$ of classes of square $-2$ as the group $\G_\D$ generated
by the reflections on hyperplanes normal to elements in $\D$.

\begin{prop}
The monodromy group $\G(\fami)$ of the smooth surface $X$ in the family
$\fami$ coincides with the group $\G_{\!_{\D_M}}$ associated to the
set
$\D_{\!M}$ of all classes of square
-$2$ contained in the sublattice $L_M\subset L_{X}$ generated by cycles in
$M=X\cap\milnor$.
\end{prop}

\proof
As in the complement of $\milnor$ the family projection is a fibre bundle by
construction, the monodromy acts only on $L_M$.
We made sure that $\milnor$ contains a Morsification for the
hypersurface singularity $z^2+y^3+x^{6\chi}$, so the set of vanishing
cycles determines generators for the monodromy group via the associated
reflections. In \cite{habil} this set is shown to coincide with the set of classes of
square $-2$ in the Milnor lattice of the singularity which is just $L_M$ and the claim
follows.
\qed

\section*{Families obtained from branched coverings}

In this section we construct families by varying the branch locus in our pull back
construction. To understand the associated vanishing cycles we consider first the
situation locally on the base, annuli covering a disc with varying branch locus:

\begin{lemma}
\labell{toptyp}
Let $Y_0$ be smooth and properly fibred over a holomorphic disc $V$ such that all
fibres are smooth elliptic curves except for a singular fibre of type $I_1$.
If $Y$ is the pull back of $Y_0$ along a double cover $d:A\to V$ such that $A$ is
an annulus and the two branch points differ from the base of the singular fibre of
$Y_0$, then
\begin{enumerate}
\item
$Y$ is properly elliptic fibred over $A$ with two singular fibres of type $I_1$,
\item
$Y$ is diffeomorphic to a trivial torus bundle over $A$ with a $2$-handle added on
each boundary with framing $-1$ along isotopic vertical curves,
\item
$H_2(Y)$ with the intersection pairing is a lattice isomorphic to the diagonal
lattice $(-2)\oplus (0)_{2g+1}$.
\end{enumerate}
\end{lemma}

\proof
The first assertion is obvious by base change properties. To get the second, we
decompose $Y_0$ as a trivial torus bundle with a $2$-handle added on a vertical
loop with framing $-1$, cf.\ \cite{kas}. The pull back of the torus bundle yields
the trivial bundle over an annulus, whereas the additional $2$-handle lifts to
handles which are added as claimed.\\
The trivial bundle supports the radical of the intersection lattice of rank
$2g+1$. The distinguished handles cancel with the same $1$-handle, hence their
cores glue to an embedded sphere of self-intersection $-2$.
\qed

\begin{lemma}
\labell{locdefo}
Let $Y_0$ be a smooth complex surface properly elliptic fibred over a holomorphic
disc $V$ with a single singular fibre of type $I_1$ at the origin,
let $$d:\afami_T\to V\times T$$ be a versal family of double covers of $U$ with
branch divisor of degree two,
and $\yfami_T$ be the pull back along $d$ of the trivial family $Y_0\times T$.\\
Then for any smooth surface $Y$ in $\yfami_T$ there are spheres $s^2_\pm$ embedded
in $Y$, such that
\begin{enumerate}
\item
$s^2_\pm$ are vanishing cycles for ordinary double point degenerations in
$\yfami_T$,
\item
$s^2_\pm$ map to simple arcs in the base, the union of which is isotopic to the
core of the annulus,
\item
the intersection of $s^2_\pm$ with a regular fibre is either empty or a vanishing
cycle for both curve degenerations,
\item
the intersection $s^2_+\cdot s^2_-$ is $\pm2$ with the sign depending on
orientations,
\item
dual classes to $s^2_\pm$ can be represented by cylinders mapping to a cocore of
the annulus.
\end{enumerate}
\end{lemma}

\proof
It suffices to give some family of double covers with the claimed properties. This
family is obtained by pull back from any versal family, hence the versal family
has the properties as well.\\
Consider first the family of double covers of the unit disc $D_1$ determined up
to Galois involution by the branch locus $B_t=\{(2t+1)/3,(2t-1)/3\}$:
\begin{eqnarray*}
\afami\hspace{4mm} & \to & \hspace{2mm}D_1\times D_1\,\ni t,t'\\
\searrow & & \swarrow\\
& D_1 & \ni t
\end{eqnarray*}
By a suitable biholomorphic equivalence we may regard $V$ as the unit disc with
the nodal curve of $Y_0$ at the origin.\\
Let $\yfami$ be the family of fibred surfaces obtained by pull back of $Y_0\times
D_1$ along $\afami\to D_1\times D_1$.
Then $\yfami$ has two ordinary double point degenerations at $t=\pm\frac12$.
The associated vanishing cycles can be represented by spheres $s^2_\pm$ in the
following way:\\
Let vanishing cells for the curve degeneration of $Y_0$ be given over the
segments $t=[0,\pm\frac12]$ with a single transversal intersection. Then the
pull back of the cells to
$\yfami$ are vanishing cells for the surface degeneration and yield spheres in
$Y$ as claimed in $i),ii),iii),iv)$. Also $v)$ is then immediate.
\qed

\begin{lemma}
\labell{tori}
If $t^2$ is a tube (cylinder) embedded into $Y_0$ which maps to a path in the base
connecting the points $s=\pm\frac13$ then its preimage $\tilde t^{\,2}$ in $Y$ is an
embedded torus which represents a class in the integer span of the classes of
$s^2_\pm$ if the boundary circles of $t^2$ are vanishing cycles for the curve
degeneration of $Y_0$.
\end{lemma}

\proof
The preimage of such a tube consists of two tubes with boundary circles identified,
a torus which represents an isotropic class in $H_2(Y)$ supported on the bundle
part, cf.\ lemma \ref{toptyp}.\\
By lemma \ref{locdefo} the sum of the classes represented by $s^2_\pm$ suitably
oriented represent an isotropic class. Since the only horizontal torus in the
bundle part which has the same intersection with all cylinders over a fixed cocore
of the base annulus is the torus with vertical loop isotopic to the vanishing
cycle of the curve degeneration, we are done.
\qed

Getting on to the global setting we define for each set of discrete invariants
$\chi,q,\mu$ a family
$\efami\cqm$ of elliptic surfaces of the given invariants:\\
\begin{description}
\item[$\chi\geq2$:]
Let $u_\cq:\cfami\to\P\times H_\cq$ be the universal family of branched covers of
$\P$ of degree $\chi$ and genus $q$ with simple branching only except for a
totally ramified branching at zero parametrized by the appropriate Hurwitz space
$H_\cq$.\\
Let $u'_\cq$ denote the restriction of $u_\cq$ to the open part $H'_\cq$
parametrizing covers with no branching in an arbitrary small neighbourhood
$U^\infty$ of infinity.\\
The family $\efami_\cq$ is then defined to be the pull back of $X_1$
along ${pr}_1\circ u'_\cq$. Since no branching is allowed on $U^\infty$ the family
is trivial over the preimage, so multiple fibres can be introduced there
simultaneously to yield $\efami\cqm$ for any given discrete invariants.
\item[$\chi=1$:]
Consider $u_{2,q}$ again, the universal family of hyperelliptic covers of genus
$q$ branching at zero and a varying set of further $2q+1$ distinct points of
$\P\setminus\{0\}$.
Denote by $u''_{2,q}$ the restriction to the part $H''$ of the corresponding
Hurwitz space which parametrizes covers with no branching in a small neighbourhood
$U^\infty$ of infinity nor in a small pointed neighbourhood $U^0_*$ of zero.
We get a family $\bar\efami_q$ as the pull back of
$\Xhalf$ along ${pr}_1\circ u'_q$.
By the triviality over the preimage of $U^\infty$ logarithmic transforms can again
be performed simultaneously.
The triviality over the preimage of $U^0=U^0_*\cup\{0\}$ on the other hand makes
it possible to normalize and blow down simultaneously and so we end up with a
family $\efami_q^\mu$ of elliptic surfaces with the preassigned invariants.
\end{description}
The choices involved in the constructions of $\efami$ and $X$ allow to assume that
the surface $X$ is contained in $\efami$, because $U^\infty$ can be chosen small
enough to exclude the branch points $x_i$ involved in the construction of $X$ and
the logarithmic transforms on $X'$ and $\efami_\cq$ can be performed in a compatible
way too.\\

We may assume that a base point $x_0$ in $\P$ is chosen such that the fibre $f_0$
in $X_1,\Xhalf$ respectively is smooth. Then there is a pair of vanishing cells
projecting to paths $p^+,p^-$ with associated vanishing
cycles in
$f_0$ dual to each other.
On the base of the elliptic surface $X_1$($\Xhalf$) we may now choose the subsets
considered above even more specifically:
\begin{enumerate}
\item
a neighbourhood $U^\infty$ of infinity not containing any singular values for the
fibration map,
\item
a simply connected neighbourhood $U$ as in prop.\ \ref{e-ref},
\item
points $x_j$,${j\geq2}$ ordered on a circle on $\P$ of constant modulus and a
corresponding geometrically distinguished system of paths $p_j$ in the complement
of $U\cup U^\infty$,
\item 
the base point $x_0$ on $\del\bar{U}$
\item
paths $p^+,p^-$ contained in $\bar{U}$
\end{enumerate}
In addition we will consider a collection $V^+_j,V^-_j,{j=2,...,2q+2}$ of
holomorphic images of a disc such that $V^+_j,V^-_j$ contain the circle
segment between $x_j$ and $x_{j+1}$, the paths $p_j,p_{j+1}$ and $p^+$, resp.\
$p^-$. Moreover they are chosen to be
disjoint from the singular values and the subset $U^\infty$.

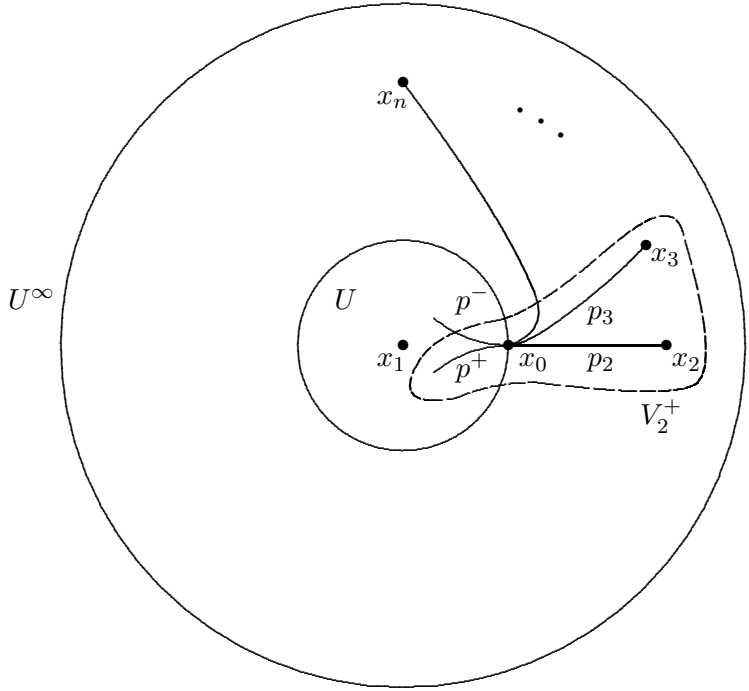
\begin{center}
\setlength{\unitlength}{7mm}
\begin{picture}(14,13)
\put(6,6){\bigcircle{4}}
\put(6,6){\bigcircle{13}}
\curvesymbol{\phantom{\circle*{.2}}\circle*{.2}}
\put(6,6){\arc[-1](5,0){22.5}}
\curvesymbol{\phantom{\circle*{.1}}\circle*{.1}}
\put(6,6){\arc[-2](3,4){10.5}}
\put(6,11){\circle*{.2}}
\put(6,6){\circle*{.2}}
\put(8,6){\circle*{.2}}
\put(9.5,5.6){$p_2$}
\put(9.5,6.5){$p_3$}
\put(8.2,5.6){${x_0}$}
\put(11.1,5.6){${x_2}$}
\put(5.5,5.6){${x_1}$}
\put(10.7,7.6){${x_3}$}
\put(5.5,10.6){${x_n}$}
\put(7,5.4){${p^+}$}
\put(7,6.7){${p^-}$}
\put(10.5,4.5){${V_2^+}$}
\put(4.7,6,7){$U$}
\put(-1.5,6.7){${U^\infty}$}
%\put(2,6.7){${\P\setminus U^\infty}$}
\put(8,8.2){\arc(0,-2.2){-40}}
\put(8,3.8){\arc(0,2.2){40}}
\put(0,0){\tagcurve(7,6.5, 8,6, 9,6.5, 10.6,7.9,  14,12)}
\put(0,0){\tagcurve(7,6.5, 8,6, 8.5,7, 6,11, 4,13)}
\put(8,6){\line(1,0){3}}
%\put(0,0){\closecurve(6.3,5.5, 6.7,5.2, 7.5,5.3, 8,5.4, 11,5.4, 11.5,7,  10.5,8.3)}
\curvesymbol{}
\curvedashes[1mm]{0,2,.5}
\put(0,0){\tagcurve(6,6, 8.5,5.3, 11.6,5.3, 11.3,8.2, 10.7,8.4, 9.5,7.5, 8,6)}
\put(0,0){\tagcurve(10,6.3, 7.9,6.5, 6.3,5.7, 7,5, 8,6)}
\put(0,0){\tagcurve(6,6.5, 7.9,6.5, 9.5,7.5,11,9)}
\put(0,0){\tagcurve(6,4.2, 7,5, 8.5,5.3, 11,5)}
\end{picture}
\end{center}

The vanishing cycles associated to these discs by lemma \ref{locdefo} then define
a set $\D^b$ of classes of square $-2$ for $X$.

\begin{lemma}
\labell{bmono}
Let $\Gcq$ be the image of the monodromy representation of $\efami\cqm$ on the
homology lattice $L^2_X$ of $X$ fibred over $C^X_q$ and let $Z^c$ be the
complement of $Z'_0$ in $X$.\\
Then with the subset $\D^b$ of $L_X^2$ of elements of square $-2$
\begin{enumerate}
\item
$\G_{\D^b}$ is a subgroup of $\Gcq$,
\item
$H_2(Z^c)$ maps into the sublattice generated by $\D^b$ and the fibre classes.
\end{enumerate}
\end{lemma}

\proof
The family $u'_\cq$ (resp.\ $u''_{2,q}$) contains a versal family of double covers
of each disc
$V^\pm_j$ and we may apply lemma \ref{locdefo} to see that all elements of $\D^b$
are in fact vanishing cycles for ordinary double point degenerations of $X$ in
$\efami=\efami^\mu_\cq$.\\
The tubes over the circle segments with the vanishing cycles as vertical
components lift thus to tori representing classes in the integer span of the
vanishing classes of the family $\efami$.
It may be checked that these tori form a basis of a unimodular lattice of rank $2q$
and thus together with the fibre classes generate $H_2(Z^c)$.
\qed

\begin{lemma}
\labell{vanlat}
Let $\D^s$ be the set of $-2$ classes of $L_X$ supported on $M=Z_0\cap\milnor$ and
$\D$ the union $\D^s\cup\D^b$ then
\begin{enumerate}
\item
$\D$ generates together with the fibre classes the orthogonal complement of the
fibre class in $L_X$,
\item
$\D$ is contained in a single $\G_\D$ orbit.
\end{enumerate}
\end{lemma}

\proof
Consider the Mayer-Vietoris sequence associated to the decomposition $X=Z^c\cup
Z'_0$. We get a map $H_2(Z^c)\oplus H_2(Z'_0)\to H_2(X)$ with a torsionfree
cokernel of rank one, which can be identified with the intersection with the class
of a general fibre. Thus the sublattice orthogonal to the fibre classes is the
image. Since $\D^b$ and the fibre classes generate $H_2(Z^c)$ and on the other
hand $\D^s$ and the fibre classes generate $H_2(Z'_0)$ claim $i)$ follows.\\
By construction each $\d\in\D^b$ intersects algebraically trivial with the tori
associated to the same vanishing cycle in $f_0$. They span a maximal isotropic
sublattice, so $\d$ is the orthogonal sum of an isotropic element in the
sublattice spanned by the tori and an element of the orthogonal complement which is
supported in $Z'_0$ and of square $-2$. Hence it differs at most by a multiple of
the fibre class from an element in $\D^s$. Hence there is another element $\d'$ in
$\D^s$ with $\d.\d'=1$. So each element of $\D^s$ is conjugated to one of $\D^s$
since $\cg_\d\circ\cg_\d'(\d)=\d'$ for $\d,\d'\in\D$ with $\d.\d'=1$. The claim
$ii)$ now holds since $\D^s$ is even a single orbit of the $\G_{\D^s}$ action.
\qed

\begin{prop}
Let $\G_X$ be the image of the monodromy representation for the surface $X$. Then
there is a set $\D'$ of homology classes of square $-2$ and a
sublattice $L'<L_X$
containing a unimodular one of corank two, such that $L'$ and the fibre classes
generate the orthogonal complement of the fibre class, $\G_{\D'}$ is a subgroup of
$\G_X$ and $L',\D'$ is a complete vanishing lattice in the sense that, cf.\
\cite[5.3.1]{habil}:
\begin{enumerate}
\item
$\D'$ generates $L'$,
\item
$\D'$ is a single $\G_{\D'}$ orbit,
\item
$\D'$ contains six elements the intersection diagram of which is\\[6mm]
\unitlength1.6cm
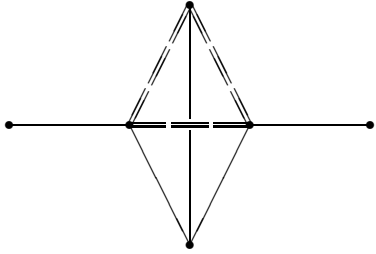
\begin{picture}(4,1.6)(-1,1.2)
\put(1,2){\circle*{.06}}
\put(1,2){\line(1,0){1}}
\put(2,2.015){\line(1,0){.3}}
\put(2.35,1.985){\line(1,0){.3}}
\put(2.35,2.015){\line(1,0){.3}}
\put(2.7,1.985){\line(1,0){.3}}
\put(2.7,2.015){\line(1,0){.3}}
\put(2,1.985){\line(1,0){.3}}
\put(3,2){\line(1,0){1}}
\put(2,2){\circle*{.06}}
\put(2.01,1.98){\line(1,2){.15}}
\put(1.98,2){\line(1,2){.15}}
\put(2.185,2.33){\line(1,2){.15}}
\put(2.155,2.35){\line(1,2){.15}}
\put(2.36,2.68){\line(1,2){.15}}
\put(2.33,2.7){\line(1,2){.15}}
\put(2,2){\line(1,-2){.5}}
\put(3,2){\circle*{.06}}
\put(2.845,2.35){\line(-1,2){.15}}
\put(2.815,2.33){\line(-1,2){.15}}
\put(2.67,2.7){\line(-1,2){.15}}
\put(2.64,2.68){\line(-1,2){.15}}
\put(3.02,2){\line(-1,2){.15}}
\put(2.99,1.98){\line(-1,2){.15}}
\put(3,2){\line(-1,-2){.5}}
\put(4,2){\circle*{.06}}
\put(2.5,3){\circle*{.06}}
\put(2.5,3){\line(0,-1){.95}}
\put(2.5,1){\circle*{.06}}
\put(2.5,1){\line(0,1){.95}}
\end{picture}
\end{enumerate}
\end{prop}

\proof
Let $\D'$ be $\G_\D\cdot\D$ then $\D'$ generates a sublattice $L'$ as claimed and
$i),ii)$ are obvious from previous results, as is $iii)$ in case $\chi\geq 2$.\\
For $\chi=1$ we note that $\D^s$ generates a sublattice of type $E_8$, hence
$\D^s$ and either half of $\D^b$ -- corresponding to the superscript '$+$' or '$-$'
-- generate a semidefinite lattice which meets the hypotheses of \cite[II.5.9]{fm}.
We may conclude that with any element of
$\D^s$ there are all elements of $L'$ in $\D'$ which differ from the given by only a
torus as considered in the proof of lemma \ref{bmono}.\\
The first four of
them $\tilde t^+_1,\tilde t^-_1,\tilde t_2^+,\tilde t_2^-$ intersect like a
symplectic basis with the only non-zero intersection being $\tilde
t_1^+.\tilde t_2^-=-\tilde t_1^-.\tilde t_2^+=1$. With elements $\a_1,\a_2,\a_3,\b$
of
$\D^s$ orthogonal except for $\a_1.\b=1$ the following six elements are in $\D'$ and
intersect as in
$iii)$:
$$
\a_2+\tilde t_1^+,\a_1+\tilde t_2^-,\a_1,\b,\a_1+\tilde t_1^-,\a_3-\tilde t_2^+.
$$
\qed

\begin{prop}
\labell{vlatt}
Let $X$ be a minimal elliptic surface with $p_g\geq q>0$. Then the set $\D$ of
classes in $L_X$ of square $-2$ orthogonal to the fibre class $f$ and represented
by a vanishing cycle in a degeneration family of $X$ generates the orthogonal
complement $L$ of the fibre class in $L_X$ and $L,\G_\D\cdot\D$ is a complete
vanishing lattice.
\end{prop}

\proof
As in \cite{ellmon} the proof relies on the deformation equivalence through
elliptic surfaces of any surface of specified discrete invariants to our reference
surface and on two further ingredients, a complete vanishing lattice as given by
the previous proposition and for each fibre class a pair of vanishing cycles the
sum of which is represented by the fibre.\\
The latter is obtained by a slight generalization to the irregular case of results
proved in \cite{ellmon}.
\qed

\section*{Monodromy results}

Monodromy defines a group of isotopy classes of orientation preserving
diffeomorphisms.
This monodromy group has natural representations on
$L^2_X:=H_2/\!\:\!$\raisebox{-0.7mm}{torsion} and on
$L^1_X:=H_1/\!\:\!$\raisebox{-0.7mm}{torsion}.\\
The image of the first is contained in the group $O(L^2_X)$ of orthogonal
transformations of $L^2_X$ with respect to the symmetric intersection pairing $q_X$.
In case of an irregular elliptic surface we have an analogous property for the
other representation:

\begin{lemma}
\labell{first}
Let $X\to C$ be an irregularly fibred minimal elliptic surface with positive
Euler number. Then the induced map on homology
$H_1(X)\to H_1(C)$ factors through an isomorphism 
$$
L^1_X\stackrel{\simeq\mbox{\hspace{10pt}}}\tto H_1(C),
$$
and the pull back of the intersection product on $H_1(C)$ is a skew-symmetric
non-degenerate bilinear form $b_X$ on $L_X^1$
such that the image of the natural monodromy representation is contained in the
group of symplectic transformations 
$$
Sp(L^1_X):=Sp(L^1_X,b_X).
$$
\end{lemma}

\proof
The fibration map is the Albanese mapping of $X$, hence induces an isomorphism
on the first rational homology. Since the integer homology of $C$ is torsionfree,
the first claim follows.\\
Furthermore the bilinear form $b_X$ is then skew-symmetric and non degenerate.\\
Finally each isotopy class of monodromy diffeomorphisms is represented by a
fibration preserving one, hence commutes with a diffeomorphism of the base.
The claim will follow as soon as the latter is shown to be orientation
preserving.\\
If it were orientation reversing, so would be the action on a
homology class of a multi\-section. On the other hand this contradicts the fact
that fibre class and intersection product are preserved under monodromy, hence the
claim.
\qed

\begin{thm}
Let $X$ be a minimal elliptic surface with $p_g\geq q>0$ and fibre class $f$. Then
the group of orthogonal transformations of its homology lattice $L^2_X$ generated
by the monodromy groups of all families containing $X$ is
$$
O'_f(L_X):=\{\cg\in O(L_X^2)|\cg(f)=f,\cg\text{ has positive real spinor norm}\}.
$$
\end{thm}

\proof
Again the argument is the same as in \cite{ellmon} provided that the ingredients
of the proof are extended to the irregular case, as we summarized in prop.\
\ref{vlatt}.
\qed

\begin{cor}
\label{triv}
Let $X$ be as above. Then $O'_f(L^2_X)$ is even generated by elements associated to
monodromy transformations acting trivially on $L^1_X$.
\end{cor}

\proof
Any diffeomorphism isotopic to the inversion at a $-2$ sphere acts trivially on
$L^1_X$, hence the claim follows from the proof of the theorem.
\qed

\begin{thm}
Let $X$ be as above. Then the group of transformations of $L^1_X$ generated by the
monodromy groups of all families containing $X$ is
$$
Sp(L_X^1).
$$
\end{thm}

\proof
Let $\cfami$ be a family of curves of genus $q$ with the full symplectic
monodromy. Then a family of elliptic surfaces deformation equivalent to $X$ is
obtained by the following construction:
Choose a divisor on $\cfami$ disjoint from the critical locus of fibre degree
$2\chi$ without horizontal component (a 2$\chi$-section).
Take modulo a suitable base change the corresponding double cover $\tilde\cfami$ and
associate the family $\tilde\cfami$ of trivial elliptic surfaces to it.
Divide out the diagonal action of the Galois action on $\tilde\cfami$ and the
involution on $E$ and
resolve the singularities of the quotient. Get thus the desired family of elliptic
surfaces over the base of $\cfami$.
Multiple fibres can again be introduced.\\
The resulting monodromy of the surface acts surely by $Sp$ on the first homology of
the base and hence on $L^1_X$ by the isomorphism of lemma \ref{first}.
\qed

\begin{cor}
Let $X$ be as above. Then the monodromy representation on $L^1_X\oplus L^2_X$ maps
onto
$$
O'_f(L_X^2)\times Sp(L_X^1).
$$ 
\end{cor}

\proof
Let $\cg$ be any element of $O'_f(L_X^2)\times Sp(L_X^1)$, then we can find
$\cg^1$ acting on $L^1_X$ as $\cg$ does. Moreover $\cg^1$ acts on $L^2_X$ as an
element of $O'_f$, hence $\cg\circ(\cg^1)\inv$ is some element of
$O'_f\times\{\id\}$ which is given by a monodromy element according to cor.\
\ref{triv}.
\qed

\begin{description}
\item[Remark:]
In \cite{kml} the image of the natural representation of the diffeomorphism group of
an elliptic surface as above on the homology lattice was shown to be generated by
$O'_f(L_X)$ and a diffeomorphism $\s$  --  induced by complex conjugation -- such
that $O'_f(L_X)$, which is now shown to coincide with the monodromy group, is a
subgroup of index two.
\end{description}

\end{document}